
\magnification = \magstephalf
\tolerance=10000              
\hbadness=10000
\font\bigtenrm=cmr10 scaled\magstep2

\font\tengo=eufm10 
\font\sevengo=eufm7 
\font\fivego=eufm5
\newfam\gothic
\textfont\gothic=\tengo 
\scriptfont\gothic=\sevengo 
\scriptscriptfont\gothic=\fivego

\def\R{{\bf R}}   
\def\Z{{\bf Z}}   
\def\Q{{\bf Q}}   
\def\C{{\bf C}}   
\def\H{{\bf H}}   

\def\Map{{\rm Map}}

\def \gn{ G_n(V^{\infty})}
\def \gr#1{ G_{#1}(V^{\infty})}
\def \g#1{G_{#1} (\H ^{\infty})}
\def \en{ E_n(V^{\infty})}

\def\sqr#1#2{{\vcenter{\vbox{\hrule height.#2pt 
         \hbox{\vrule width.#2pt height#1pt \kern#1pt 
                 \vrule width.#2pt} 
         \hrule height.#2pt}}}}

\def\A{{\cal A}} 
\def\AH{{\cal A}_H}
  
\def\B{{\cal B}} 
 
\def\BH{{\cal B}_H}
\def\BHO{{{\cal B}_H^0}} 

\def\BO{{\cal B}^0}
\def\G{{\cal G}}

\def\GH{{{\cal G}_H}} 
\def\GHO{{{\cal G}_H^0}} 
\def\GO{{{\cal G}^0}}
\def\M{{\cal M}} 
 
\def\MH{{\cal M}_H}

\def\Fg{{\cal F}_{2g}} 
\def\FS{{\cal F}_{2g} \times_{\sigma} S(v) } 
\def \Mg{M_{2g}} 


\centerline {\bigtenrm Morse theory for the Yang-Mills functional}
\smallskip
\centerline {\bigtenrm via equivariant homotopy theory}

\smallskip 

\centerline { U. Gritsch}

\smallskip 



\proclaim Abstract. 

In this paper we show the existence of non minimal critical points of
the Yang-Mills functional over a certain family of 4-manifolds $\{
M_{2g} : g=0,1,2, \ldots \}$ with generic $SU(2)$-invariant
metrics using Morse and homotopy theoretic methods.  These manifolds
are acted on fixed point freely by the Lie group $SU(2)$ with quotient
a compact Riemann surface of even genus. We use a version of invariant
Morse theory for the Yang-Mills functional used by Parker in [Pa1] and
R\aa de in [R\aa].

\proclaim 1~Introduction. 

In this paper we construct a family of Riemannian spin {4-\nobreak
manifolds}, denoted by $\{ M_{2g} : g=0,1,2, \ldots \}$, acted on by
the Lie group $SU(2)$ and prove an existence result for non minimal
critical points of the Yang-Mills functional over the manifold $\Mg$
(for $g \ne 1$ ) for generic $SU(2)$-invariant metrics.  
The manifold $\Mg$ is the product of a compact Riemann surface $\Fg$ of 
even genus and the two-sphere $S^2$ modulo an involution. The $SU(2)$-action 
is the standard action on $S^2$ and the trivial action on $\Fg$. 
In section 3.4 we prove

\proclaim Theorem 3.4.10. 
Fix $g=0,1,2,\dots$ and let $\Delta^+$ and $\Delta^-$ be 
the positive and negative $SU(2)$-equivariant spinor bundles over the 
4-manifold $\Mg$.  
Fix an $SU(2)$-invariant metric. 
\item{(i)} The Yang-Mills functional restricted to the invariant orbit space 
$ {\cal YM}: \B_{SU(2)}\rightarrow \R$  
over the manifold $\Mg$ has at least $2g+1$ critical points on each of 
the bundles $\Delta^+$ and $\Delta^-$. 
\item{(ii)} In the case $g=0$ for a generic $SU(2)$-invariant 
metric the 
critical point on the bundle $\Delta^+$ cannot be self dual and on the bundle 
$\Delta^-$ it cannot be anti self dual. 
In the case $g\geq 2$ for a generic $SU(2)$-invariant metric 
at least one of 
the critical points on the bundle $\Delta^+$ cannot be anti self dual and at 
least one of the critical points on the bundle $\Delta^-$ cannot be self dual. 

By Palais' principle of symmetric criticality [Pal2] 
each of the critical points in (i) is a critical 
point of the Yang-Mills functional on the non-equivariant orbit space  
${\cal YM}:\B \rightarrow \R $,  
i.e. they are Yang-Mills connections. 

For a long time it was conjectured that by analogy with the harmonic
map problem on the two-sphere the Yang-Mills functional on the
4-sphere does not have other critical points besides the self dual and
anti self dual connections.  In [SSU]
L.M. Sibner, R.J. Sibner and K.K. Uhlenbeck however construct a
sequence of non minimal critical points on the trivial $SU(2)$-bundle
on the 4-sphere with the round metric. In [SS] L. Sadun and J. Segert
find a critical point of the Yang-Mills functional on every
$SU(2)$-principal bundle $\eta$ over the 4-sphere $S^4$ with the
standard metric provided its second Chern class $c_2 (\eta )$ is not
equal to $\pm 1$.  In [Wan] H.-Y. Wang proves the existence of an
infinite number of irreducible $SU(2)$-connections over the manifolds
$S^2 \times S^2$ and $S^1 \times S^3$ with the standard metrics which
are non minimal solutions to the Yang-Mills equations.  Finally in
[Pa2] T.H. Parker proves the existence of irreducible non-minimal
Yang-Mills fields on the trivial bundle over $S^1_L \times S^3$ for
some $L$, where $S^1_L$ denotes the circle of radius $L$; and he shows
that there exists a family of metrics on the 4-sphere $S^4$ each of
which admits an irreducible Yang-Mills field on the trivial
$SU(2)$-bundle.  

Naturally one would like to do Morse theory for the Yang-Mills
functional defined on the space of connections on a principle bundle over a 
Riemannian 4-manifold 
modulo the action of the gauge group.  Since this space $\B$ is
infinite dimensional the Yang-Mills functional needs to satisfy a
certain compactness condition, called the Palais-Smale
condition. However, by Uhlenbeck's weak compactness theorem, it is
well-known that the Yang-Mills functional defined on the orbit space
$\B$ does not satisfy this condition. 

As noted in [Pa1] this problem does not occur if the 4-manifold $M$ is
acted on by a compact Lie group $H$ such that the action has no zero
dimensional orbits and we restrict the Yang-Mills functional to the
orbit space $\BH$ of $H$-invariant connections modulo the
$H$-invariant gauge group.  

A proof was outlined in [Pa1] that the Yang-Mills functional defined
on the invariant orbit space $\BH$ satisfies the Palais-Smale
condition provided the orbit space $\BH$ has no singularities.
Unfortunately this proof contains a gap: It is not clear that one can
choose good invariant gauges. A different proof was given recently by
J. R\aa de in [R\aa ]. (See also section 3.1 of this paper). 

In [Pa1] Parker uses his Morse theory to reprove Sadun's and Segert's
result in [SS] using a fixed point free $SU(2)$-action on $S^4$ called
the quadrupole action.  There is also a fixed point free $SU(2)$-action on
$S^2 \times S^2$ but Parker's Morse theory is too weak to establish
the existence of irreducible non (anti) self dual Yang-Mills
connections. 
The manifold $\Mg$ in theorem 3.4.10 is acted on fixed point freely 
by the Lie group $SU(2)$, too. 

The general strategy 
to prove theorem 3.4.10 
as proposed by Parker in [Pa1] 
is as follows. 
Let $\MH \hookrightarrow \BH$ denote the subspace of the 
absolute minima of the Yang-Mills functional modulo the 
invariant gauge group. We assume that the second Chern class of 
the bundle we are working with is positive. Then the absolute 
minima of the Yang-Mills functional are the anti self dual 
connections and there are no self dual connections. 
(By changing the orientation we obtain the opposite case of 
negative second Chern class). 
If the inclusion $\iota : 
\MH \hookrightarrow \BH$ is not a homotopy equivalence 
then the Palais-Smale condition for the Yang-Mills 
functional implies the existence of at least one non minimal critical 
point of the Yang-Mills functional on the space $\BH$. For otherwise the retraction 
along the flow lines 
of the gradient flow of the Yang-Mills functional 
would give a homotopy inverse of the inclusion 
$\iota : \MH \hookrightarrow \BH$. 

It follows from Uhlenbeck's weak convergence theorem and its
equivariant versions in [Cho] and [Ba] that if the bundle $\eta$ does
not admit a reducible $H$-invariant anti self dual connection for any
metric then for a generic $H$-invariant metric 
the moduli space $\MH$ is a (possibly
empty) compact closed manifold of a dimension which can be computed
using the equivariant Atiyah-Singer index theorem.

Denote by $\BHO$ the based invariant orbit space of invariant
connections modulo the based invariant gauge group.  In section 2 of
this paper we identify the weak homotopy type of the space $\BHO$ in
terms of an invariant mapping space. This is a direct generalization
from the corresponding non equivariant theorem in [AB, prop. 2.4,
p.540] or [DK, prop. 5.1.4, p.174].  It will turn out that in our
example of theorem 3.4.10 (for $M=\Mg$ and $H=SU(2)$) the orbit spaces
$\BH$ and $\BHO$ are homeomorphic. Then the information on the
homotopy type of the space $\BH$ and the dimension of the manifold
$\MH$ for generic $H$-invariant metrics will be enough to deduce that
the inclusion $\MH \hookrightarrow
\BH$ cannot be a homotopy equivalence.

In section 3.2 we define the manifolds $\Mg$ for $g=0,1,2, \ldots $ acted 
on fixed point freely 
by the Lie group $SU(2)$, show that they are spin and compute the second 
Chern classes of the positive and negative spinor bundles 
$\Delta^+$ and $\Delta^-$.
In section 3.3 we compute the homotopy type of the invariant orbit space 
$\B_{SU(2)}$ and in section 
3.4 we prove theorem 3.4.10. 

\noindent{\sl Acknowledgments:} I would like to thank 
my advisor Ralph Cohen for suggesting the problem and for constant
support and encouragement.  I thank the Studienstiftung des deutschen
Volkes for a dissertation fellowship.  This paper is part of 
my Ph.D. thesis written at Stanford University, 1997.

\bigskip  

\proclaim 2~The homotopy type of the space $\BHO$.

\proclaim 2.1 Notation and terminology. 

Let $M$ be a closed oriented Riemannian 4-manifold and let 
$ \eta = (P \rightarrow M)$ be an $SU(2)$-principal bundle 
over the manifold $M$. Denote by $\A$ the space of 
smooth connections on the bundle $\eta$. Let $\G$  be the group of 
smooth gauge transformations of this bundle. This means that elements 
$g \in \G$ are smooth automorphisms of the bundle $\eta$. Let 
$\GO$ be the subgroup of the group $\G$ of gauge transformations whose 
elements are 
the identity over a given base point $m \in M$. 

Let $H$ be a compact Lie group acting smoothly on the manifold 
$M$ on the left. We also assume that the action lifts to the 
total space $P$ of the bundle $\eta$ 
such that the left action of the Lie group $H$ on the total 
space $P$ commutes with the right action of the structure group 
$SU(2)$ on $P$. 
We define the invariant gauge group $\GH$  
to be the subgroup of the gauge group $\G$ such that each gauge 
transformation commutes with the action of the Lie group $H$ 
on the total space $P$ of the bundle $\eta$. 
Analogously we define the based invariant gauge group 
$\GHO$ to be the equivariant gauge transformations which are the 
identity over a given based orbit $ {\cal O}_m \subset M$. 

The groups $\G$ and $\GO$ act naturally on the space of connections 
$\A$ from the right by the pull back of connections. 
We define the orbit space 
$ \B $ to be the space $ \A / \G $ of connections modulo the action of the 
gauge group and the based orbit space $\BO $ to be the space $ \A / \GO $  
of connections modulo the action of the based gauge group. 

The left action of the group $H$ on the bundle $\eta$ induces 
a right action on the space of connections $\A$ again by the pull back 
of connections. The fixed points of this action are called 
$H$-invariant connections. We denote the space of fixed points 
by $\AH$. The invariant groups $\GH$ and $\GHO$ act on the space 
of invariant connections $\AH$ as above. As in the non-equivariant 
setting we define the invariant orbit space 
$ \BH $ to be the space $ \AH / \GH $ of invariant connections modulo the 
invariant gauge group and the based invariant orbit space 
$ \BHO $ to be the space $ \AH / \GHO $ of invariant connections 
modulo the based invariant gauge group.

It is customary in gauge theory to complete the space $\A$ in the $L^{2,2}$-Sobolev norm
and the gauge group $\G$ in the $L^{3,2}$-Sobolev norm 
using a fixed connection $A$ on the bundle $\eta$.  For notational
convenience we still denote the completed space of connections by $\A$
and the completed gauge group by $\G$. Then the group $\G$ is a
Hilbert Lie group which acts smoothly on the space of connections
$\A$.  Analogously we complete the space of invariant connections
$\AH$ in the $L^{2,2}$-Sobolev norm and the invariant gauge group $\GH$
in the $L^{3,2}$-norm using an $H$-invariant connection which we fix once
and for all.  The based gauge group $\GO$ is a closed subgroup (in the
$L^{3,2}$-norm) of the full gauge group $\G$. The same holds for the
based invariant gauge group $\GHO$. The orbit spaces $\B$ and $\BO$
and the invariant orbit spaces $\BH$ and $\BHO$ inherit a topology
from the $L^{2,2}$-norm on the space of connections $\A$ or the space of
invariant connections $\AH$.

The results of section 2 are also true if we complete the spaces $\A$
and $\AH$ using Morrey norms as used in [R\aa ].  This will be used in
section 3 and explained in 3.1.  The reader can also take the spaces
$\B$ and $\BH$ to be the smooth (invariant) connections modulo the
smooth (invariant) gauge group together with a topology induced from
Sobolev or Morrey norms on the space of smooth (invariant)
connections.

\proclaim 2.2 Classification of equivariant bundles. 

In order to get a model for the space $ \BHO $ in terms of an 
equivariant mapping space we need to classify 
equivariant principal bundles. These 
are classified in much the same way as ordinary bundles are. 
For computational reasons which will become clear later 
we prefer to describe the classification theory of equivariant 
vector bundles. However 
the classification theories of equivariant vector bundles 
and equivariant principal bundles are really equivalent. 
We only give definitions and state results. 

Let $V$ be a fixed unitary $H$-representation. 
Let $ G_n(V^k)$ be the Grassmannian of unitary 
$n$-planes in the vector space 
$V^k = \underbrace{V \oplus \dots \oplus V}_{k \ \rm times} $ 
for $k >n$. The group $H$ acts naturally on the space 
$ G_n(V^k)$ because it acts on the space $V^k$.  
We denote by 
$\gn = \lim_{k \rightarrow \infty } G_n(V^k)$ the direct limit as 
an $H$-equivariant space.  
Let $ \gamma_n (V^k ) = ( E_n (V^k) \rightarrow G_n (V^k ) )$ 
denote the canonical unitary bundle over the 
Grassmannian. The fiber over a plane $P \in G_n (V^k )$ 
are the points $ p \in P$. 
The $H$-action on the Grassmannian lifts naturally to an $H$-action on the 
total space of the canonical bundle giving this bundle the structure of an 
$H$-equivariant bundle. We take the limit 
$ \en = \lim_{k \rightarrow \infty } E_n(V^k)$
and obtain the $H$-equivariant bundle 
$ \gamma_n (V^{\infty}) = ( \pi : \en \rightarrow \gn )$. 

Define $ {\rm Vect}^{V,n}_H (M) $ to be the isomorphism 
classes of 
$H$-equivariant unitary $n$-dimensional vector bundles over $M$ 
with the following property: 
For every $m \in M$ the isotropy representation 
of the isotropy group $H_m$ on the fiber over $m$ is equivalent to 
a sub-module of the  
$H_m$-module $V^n$ induced by the given 
$H$-module $V$. We call such vector bundles ``subordinate 
to the representation $V$''. 
Let $ [M, \gn ]^H $ denote the $H$-equivariant homotopy classes 
of the $H$-maps from the manifold $M$ to the equivariant 
Grassmannian $\gn$. 

\proclaim Theorem 2.2.1. [Wasserman, Segal]
The map 
$$ 
\eqalign{ [M, \gn ]^H & \rightarrow {\rm Vect}_H^{V,n} (M) \cr 
          [f] & \mapsto (f^{\ast}(\gamma_n (V^{\infty}))) \cr } 
$$ 
is well defined, natural in $M$ and an isomorphism. The same 
classification theorem holds in the case of orthogonal or symplectic 
equivariant vector bundles if we choose the $H$-module $V$ to be 
orthogonal or symplectic. 

\noindent {\it Proof.} A proof can be found in [Wa, section 2, p.132] in the case 
of orthogonal equivariant vector bundles. The same proof carries over 
to unitary and symplectic equivariant vector bundles. A proof of a similar 
result in the case of complex vector bundles 
over a compact manifold is also given in 
[Se, section 1, p.131].

We need a similar classification theorem for bundles with a fixed
equivariant trivialization over a fixed orbit.  Fix one orbit ${\cal
O}_m \cong H/H_m$ for some $m \in M$ once and for all. Let $ {\rm
Vect}_{H,0}^{V,n}(M) $ denote the set of isomorphism classes of
equivariant $n$-dimensional unitary vector bundles over $M$
subordinate to the representation $V$ together with a fixed
trivialization over the orbit ${\cal O}_m$. Fix a point $\ast \in
\gn$.  Let $[M, \gn]^H_0$ denote the equivariant pointed homotopy
classes of $H$-maps $f: M \rightarrow \gn $ which map the point $m$ to
the chosen point $ \ast \in \gn$ and hence which map the orbit ${\cal
O}_m \subset M $ to the orbit $ {\cal O}_{\ast} \subset \gn$.

\proclaim Theorem 2.2.2. 
Assume that the manifold $M$ is compact. 
There is a point $\ast \in \gn$ such that 
the map
$$ 
\eqalign{ [M, \gn ]^H_0&\rightarrow  {\rm Vect}_{H,0}^{V,n} \cr 
          [f]& \mapsto ( f^{\ast} (\gamma_n (V^{\infty})))\cr } 
$$ 
is well defined, natural in $M$ and an isomorphism. 
The same 
classification theorem holds in the case of orthogonal or symplectic 
equivariant vector bundles if we choose the $H$-module $V$ to be 
orthogonal or symplectic. 

\noindent {\it Proof.} 
The proof of theorem 2.2.2 uses the same methods as the proof of theorem 2.2.1
working in the category of based equivariant vector bundles.  

\bigskip  

\proclaim 2.3 The weak homotopy equivalence $\BHO \simeq 
\Map ^0_H (M, \, B(H, Sp(1)))^{\eta}$.  

\medskip 

Recall that $ \eta = (P \rightarrow M)$ is an $H$-equivariant 
$SU(2) \cong Sp(1)$-principal bundle and we also denote by $\eta = (E= P \times_{Sp(1)} \H 
\rightarrow M)$ the associated quaternionic line bundle with 
structure group $SU(2) \cong Sp(1)$. 
Choose an orbit ${\cal O}_m \cong H/H_m$ in $M$ and a fixed 
trivialization of the bundle $\eta$ over ${\cal O}_m$. 
Choose a quaternionic representation $V$ of the group $Sp(1)$ 
such that the vector bundle $\eta$ is subordinate to $V$. 
Denote by $B(H,Sp(1))= \gr {1} $ the corresponding 
Grassmannian which classifies the vector bundle $\eta$ together with 
the trivialization over the orbit ${\cal O}_m$ according to 
theorem 2.2.2. 
In the following if an $H$-equivariant vector bundle $\xi$ is 
subordinate to an $H$-module $V$ 
we say that the associated principle frame bundle (also denoted by $\xi$) 
is subordinate to the $H$-module $V$ as well. 
Also, throughout this paper, we give all spaces of maps the compact-open 
topology. 

To prove the weak homotopy equivalence 
$$
\BHO \simeq 
\Map _H^0 (M , \; B(H,Sp(1)))^{\eta} 
$$ 
(where the right hand side denotes the component of the mapping space of maps 
which classify the bundle $\eta$) 
we follow the proof of the corresponding non-equivariant result 
$$
\BO \simeq 
\Map ^0 (M , \; B(Sp(1)))^{\eta} 
$$ 
in [DK, prop. 514, p.174]. 

As in [DK, (5.1.5), p.175] associated to the 
$H$-equivariant $SU(2)$-principal bundle 
$\eta = (P \rightarrow M)$ we have the bundle $\Xi$ over the base 
$\BHO \times M$ defined as 
$ SU(2)\rightarrow \AH \times_{\GHO} P  
\rightarrow \BHO \times M $. 
The bundle $\Xi$ is in a natural way an $H$-equivariant bundle and hence 
it is classified by an $H$-equivariant map 
$
\delta : \BHO \times M \rightarrow G_1 (V^{\infty}) = B(H, Sp(1)) $. 
Here $V$ is an $H$-module to which the bundle $\eta$ (and therefore also 
the bundle $\Xi$) is subordinate. 

\proclaim Theorem 2.3.1. The adjoint of the map $\delta$ 
$$ 
\delta^{ad} : \BHO \rightarrow \Map^0_H( M , \, B(H,Sp(1)) )^{\eta}  
$$ 
is well-defined and 
induces a weak homotopy equivalence. 

\noindent {\it Proof.} The proof is a technical modification of the analogous 
result in the non-equivariant setting given in [DK, prop. 5.1.4, p. 174]. 
One uses a universal family of $H$-invariant framed connections on the 
bundle $\Xi$. 

\bigskip  

\proclaim 3~Non minimal critical points of the Yang-Mills functional. 

In the remaining part of this paper we construct the family $\{ M_{2g}
: g=0,1,2, \ldots \}$ of 4-manifolds and prove theorem 3.4.10.

\proclaim 3.1 Analytical background. 

We introduce certain completions of the spaces $\A$, $\AH$,  
$\G$ and $\GH$ 
not usually used in gauge theory but used in [R\aa , \S 4]. 
Using these norms J. R\aa de is able to show ([R\aa, Corollary 3, p. 3]) 
that the Yang-Mills functional on the space
$\BH = \AH / \GH $ 
satisfies the Palais-Smale condition provided the group $H$ acts
isometrically with no zero dimensional orbits on the manifold $M$ and
the $H$-equivariant 
principal bundle $\eta = (SU(2) \rightarrow P \rightarrow M)$ has
no reducible connections.
(A connection is called reducible if its isotropy group of the action of the 
gauge group on the space of connections is larger than the 
subgroup $\Z_2 = \{ \pm 1 \}$ of constant gauge transformations). 

The Morrey space $L^p_{\lambda} (\R^n )= L^{0,p}_{\lambda} (\R^n )$, with 
$p \in [1, \infty )$ and $\lambda \in \R$, is defined as the space of all 
$f \in L^p (\R^n )$ such that 
$$
\sup_{\rho \in (0,1]} \sup_{x \in \R^n} \rho^{\lambda -n}
\| f \| ^p_{L^p (B_{\rho} (x))} < \infty \; . 
$$ 

It is a Banach space with norm 
$$
\| f \|^p_{L^p_{\lambda}(\R^n)} = \| f \|^p_{L^p (\R^n)} + 
\sup_{\rho \in (0,1]} \sup_{x \in \R^n} \rho^{\lambda -m} 
\| f \| ^p_{L^p (B_{\rho} (x))} \; \; . 
$$ 
The Morrey space $L^{k,p}_{\lambda} (\R^n )$ with $k$ a positive integer, $p \in [1, \infty )$ 
and $\lambda \in \R$, is defined as the space of all $f \in L^{k,p}_{\lambda} (\R^n )$ 
such that 
$\partial^{\alpha} f \in L^p_{\lambda} (\R^n )$ for all $\alpha$ with 
$\vert \alpha \vert \leq k$. It is a Banach space with norm 
$$
\| f \|^p_{L^{k,p}_{\lambda}(\R^n)} = 
\sum_{\vert \alpha \vert \leq k} \| \partial^{\alpha} f \|^p_{ L^p_{\lambda} (\R^n )} \; \; . 
$$ 
The global Morrey spaces $L^{k,p}_{\lambda}(M) $ are defined using a local 
trivialization on the manifold $M$ ([R\aa , p.10]). 

The reason, why Morrey spaces are useful in equivariant gauge theory
are the following two observations made in [R\aa ].

\proclaim Lemma 3.4.1. (Special case of [R\aa, Lemma 4.1, p.11 ]) 
Let $H$ be a compact Lie group that acts smoothly on the manifold $M$, in 
such a way that all $H$-orbits have dimension $\geq 1$, and that acts smoothly 
on the bundle $\eta$. If $s \in L^{k,p} (M, \eta)$ with $k \in \Z$ and 
$p \in [1, \infty)$, is $H$-invariant, then 
$s \in L^{k,p}_s (M, \eta )$. If $A$ is a $H$-invariant connection on $\eta$, 
then 
$$ 
\| s \|^p_{L^{k,p}_{3,A}(M)} \leq c \| s \|^p_{L^{k,p}_{A}(M)} \; . 
$$
The constant $c$ only depends on the Riemannian manifold $M$ and on
the orbits of the action of the group $H$ on the manifold $M$.

The second point is that the Morrey spaces $L^{k,p}_d $ in $n$ dimensions 
satisfy ``the same'' embedding theorems as the Sobolev spaces 
$L^{k,p}$ in $d$ dimensions. Hence we get multiplications 
$$ 
\eqalign{ 
L^{2,2}_3 \times L^{2,2}_3 &\rightarrow L^{2,2}_3 \cr 
L^{2,2}_3 \times L^{1,2}_3 &\rightarrow L^{1,2}_3 \cr  
{\rm and} \;\;\; 
L^{1,2}_3 \times L^{1,2}_3 &\rightarrow L^3_3 \rightarrow L^2 \; \; \; 
\; \; \; ( \hbox {see [R\aa , p.12]}) \; . \cr  }
$$ 
It follows that if we complete the gauge groups $\G$ and $\GH$ and the spaces $\A$
and $\AH$ in the $L^{2,2}_3 $- and $L^{1,2}_3$-norm respectively then $\G$ and $\GH$ are 
Hilbert Lie groups that acts smoothly on $\A$ or $\AH$ and the Yang-Mills
functional is continuous on $\A$ and $\AH$. 

J. R\aa de is now able to show in 
[R\aa , Corollary 3, p.3] that the Yang-Mills functional on the space
$\BH = \AH / \GH $ 
satisfies the Palais-Smale condition provided the group $H$ acts
isometrically with no zero dimensional orbits on the manifold $M$ and
the $H$-equivariant 
principal bundle $\eta = (SU(2) \rightarrow P \rightarrow M)$ has
no reducible connections.

If the bundle $\eta$ does not admit reducible connections the invariant orbit 
space $\BH$ is an infinite dimensional Hilbert manifold. The $L^{1,2}$-metric on the 
space of invariant connections $\AH$ is $\GH$-invariant and descends to 
a metric on the manifold $\BH$. 
Lemma 3.1.1 shows that the space $\BH$ together with the $L_3^{1,2}$-topology and the 
$L^{1,2}$-metric is a complete Riemannian manifold.

\proclaim 3.2 The manifold $\Mg = \FS$.  

Let $\Fg$ denote an oriented Riemann surface of genus $2g$.
We denote by $S(v)$ the unit sphere in the 
representation $v$ where $v$ is the standard representation of the Lie group 
$SO(3)$ on $\R^3$. The space $S(v)$ has an induced $SO(3)$-action. 
Topologically the represenation sphere $S(v)$ is just the 2-sphere $S^2$.  
We define an involution $\sigma$ on the product $\Fg \times S(v)$ 
as follows. We think of the Riemann surface $\Fg$ as being obtained from 
two oriented Riemann surfaces ${\cal F}_g$ of genus $g$ with a 
2-disk $D^2$ removed and glued along the boundary $S^1$, i.e. 
$$ 
\Fg = ( {\cal F}_g - D^2 ) \cup_{S^1} ({\cal F}_g - D^2 ) \, . 
$$ 
We say that one of the surfaces $ {\cal F}_g - D^2 $ is positive and 
think of it as lying above the circle $S^1$. We call the other negative, 
think of it as lying below the circle $S^1$ and as the mirror image 
of the positive one. 
Then there is an involution 
$$
\sigma_1 : \Fg \rightarrow \Fg 
\eqno{(3.2.1)} 
$$ 
fixing 
the ``gluing circle'' $S^1$ and interchanging the positive and negative 
surfaces.  

\noindent {\it Example 3.2.2.} In the case $g=0$ the involution $\sigma_1$ is the 
map 
$$ 
\eqalign{ \sigma_1 : S^2 & \rightarrow S^2 \cr 
          \sigma_1 (x,y,z) &= (x,y,-z) \, . \cr }
$$ 

Let $\sigma_2 : S(v) \rightarrow S(v)$ be the reflection through 
the origin, i.e.,  
$\sigma_2 (x,y,z) = (-x,-y,-z)$. 
Then we define 
$$ 
\eqalign{\sigma : \Fg \times S(v) & \rightarrow \Fg \times S(v) \cr 
      {\rm by} \qquad   \sigma (v,w)&= (\sigma_1 (v), 
\sigma_2 (w)) \, \; \; . \cr }         
\eqno{(3.2.3)}
$$ 
Since the involutions $\sigma_1 $ and $\sigma_2$ are both orientation 
reversing the involution $\sigma : \Fg \times S(v) \rightarrow \Fg \times S(v)$ 
is orientation preserving. Also since the involution 
 $\sigma : \Fg \times S(v) \rightarrow \Fg \times S(v)$ has no fixed points 
the quotient 
$$ 
\Mg = \FS 
\eqno{(3.2.4)} 
$$ 
is a closed oriented 4-manifold. 
Any pair of metrics on the two-sphere $S(v)$ invariant under the involution 
$\sigma_2$ and on the Riemann surface $\Fg$ invariant under the involution 
$\sigma_1$ induce a metric on the manifold $\Mg$. We fix one arbitrary metric 
on the manifold $\Mg$. 

The action of the Lie group $SO(3)$ on the sphere $S(v)$ 
commutes with the reflection through the origin in $S(v)$ and hence we obtain 
an $SO(3)$-action on the manifold $\Mg$ by letting the group $SO(3)$ act 
trivially on the surface $\Fg$ and by the action on $S(v)$. 
For technical reasons which will become clear later we consider the induced 
$SU(2)$-action on the manifold $\Mg$. This action has 
no fixed points as $SO(3)$ acts 
transitively on $S(v)$. Also since the group $SU(2)$ 
is connected it acts in an orientation 
preserving fashion on the manifold $\Mg$. 
Since the group $SU(2)$ is compact we can assume that the chosen metric 
on the manifold $\Mg$ is $SU(2)$-invariant. 

\proclaim Lemma 3.2.5. The manifold $\Mg = \FS$ is a spin manifold. 

\noindent {\it Proof.} Let $T \Mg$ denote the tangent bundle of the manifold 
$\Mg$. We have to show that the second Stiefel-Whitney class
$$
\omega_2 (T\Mg ) \in H^2 (\Mg ; \Z_2 ) = 
{\rm Hom} \, (H_2 (\Mg ; \Z_2 ); \Z_2 ) 
$$ 
is zero. 
Let $\{ x_i : i\in I \} \subset H_2 (\Mg ; \Z_2 )$ be a generating set 
of $ H_2 (\Mg ; \Z_2 )$ and represent every element $x_i$ by a map 
$f_i : X \rightarrow \Mg$ where $X$ is a compact surface. It is enough to 
show that for every map $f_i$ as above 
the pull back bundle 
$f_i^{\ast} (T\Mg )$ is trivial. 
We now only sketch the argument. 

The group 
$H_2 (\Mg ; \Z_2 ) \cong 
H_2 ( \Mg ; \Z ) \otimes \Z_2 
\oplus {\rm Tor}\, (H_1 (\Mg ; \Z ) ; \Z_2 )  
         \cong \Z_2 \oplus \Z_2 $ 
is generated by two cycles 
$f_1 e_1^+ $ and $f_0 e_2^+$ chosen as follows. 
On the two-sphere $S(v)$ we choose a CW-decomposition symmetric under 
the reflection through the origin consisting of 
two 0-cells $e_0^+$ and $e_0^-$, two 1-cells 
$e_1^+$ and $e_1^-$ and two 2-cells $e_2^+$ and $e_2^-$. 
The cell $f_0$ denotes a 0-cell on the Riemann surface $\Fg$ which lies 
on the ``gluing circle''. 
Then $f_1 e_1^+ $ and $f_0 e_2^+$ denote 
the images of the cartesian 
products of the appropriate cells on the Riemann surface $\Fg$ and the 
two sphere $S(v)$ under the projection 
$\pi : \Fg \times S(v) \rightarrow \Fg \times_{\sigma} S(v) = \Mg$. 

The map 
$f: \R P^2  \rightarrow \FS =\Mg $ given by 
$ [x] \mapsto [f_0 , x] $ represents the torsion element $f_0 e_2^+$. 
One shows that $f^{\ast} (T\Mg )$ is the trivial bundle 
using the same ideas as in [MS, Lemma 3.4, p.43] to compute the 
tangent bundle of the manifold $\R P^2$. 

Let $g : S^1 \times \R P^1 \rightarrow \FS$ be the map 
$g (\lambda , [\mu ]) = [\lambda , \mu ]$ where 
we identify the circle $S^1$ with the ``gluing circle''
$S^1 \subset \Fg$. The map $g$ is a representative  
for the element $f_1 e_1^+ $. 
One shows that 
$ g^{\ast} (T\Mg ) \cong S^1 \times f^{\ast} (T \Mg ) / \R P^1 $  
where we consider the 1-dimensional real projective space 
$\R P^1$ sitting inside $\R P^2$ by the standard inclusion. 
This finishes the proof of lemma 3.2.5.

Fix a spin structure corresponding to an element $\sigma \in H^1 (\Mg ; \Z )$. 
Let $\Delta^+_{\sigma}$ and $\Delta^-_{\sigma}$ be the associated positive and 
negative spinor bundles. Both are 2-dimensional complex vector bundles 
with structure group $SU(2)$. 

\noindent {\it Remark 3.2.6.} As $SU(2)$-bundles the bundles 
$\Delta^+_{\sigma} $ and $\Delta^-_{\sigma}$ do not depend on the 
choice of the spin structure on the manifold $\Mg$. This means that for 
two elements $\sigma$ and $\delta \in H^1 (\Mg ; \Z_2 )$ if 
$\Delta^+_{\sigma} \, , \, \Delta^-_{\sigma}$ and 
$\Delta^+_{\delta} \, , \, \Delta^-_{\delta}$ 
denote the associated positive and negative spin bundles then as 
$SU(2)$-bundles
$\Delta^+_{\sigma}\cong \Delta^+_{\delta}$ and 
$\Delta^-_{\sigma}\cong \Delta^-_{\delta}$.
This can be proved using an idea which is stated in [LM, p. 84]. 

Recall that the Lie group 
$H = SU(2)$ acts on the manifold $\Mg$ by acting trivially on the surface 
$\Fg$ and by the action of $SO(3)$ on $S(v)$. 
This action on the manifold $\Mg$ has no fixed points 
and lifts canonically to the tangent bundle 
$T\Mg \cong T\Fg \times_{\sigma} TS(v) $. Since the Lie group $SU(2)$ 
is connected and simply connected it lifts to any spin bundle over the 
manifold $\Mg$ covering the action on the tangent bundle. Hence it also lifts 
to the positive and negative spinor bundles $\Delta^+$ and $\Delta^-$. 
These $SU(2)$-actions on the bundles  $\Delta^+$ and $\Delta^-$ are 
lifts of the induced $SU(2)$-actions on the bundles 
$E_+ = T\Mg \times_{\rho_+} \Lambda_+^2 (\R^4 )$ and 
$E_- = T\Mg \times_{\rho_-} \Lambda_-^2 (\R^4 )$ to a spin structure 
on $E_+$ and $E_-$. 
Here $\rho_{\pm}$ denote the two non equivalent irreducible three dimensional 
representations of $SO(4)$ on $\Lambda^2_{\pm}(\R^4)$. 
These lifts from the bundles $E_+$ and $E_-$ to the 
bundles $\Delta^+$ and $\Delta^-$ are unique since each two lifts differ by 
maps $\alpha_{\pm} : SU(2) \times P_{ Spin } (E_{\pm}) 
\rightarrow \Z_2$ such that 
 $\alpha_{\pm} ({\rm Id}, p)=1$ for each $p \in  P_{{\rm Spin }} 
 (E_{\pm})$. 
Here we denote by $P_{Spin}(E_+ )$ the principal spin bundle 
of the bundle $E_+$ and define $P_{Spin}(E_- )$ similarly. 
Since the product 
$SU(2) \times P_{Spin} (E_{\pm})$ is connected the maps  
$\alpha_{\pm}$ have to be the constant maps  
$\alpha_{\pm}(\lambda , p) = 1 \in \Z_2$ for all 
elements $\lambda \in SU(2)$ and $p \in P_{ Spin}(E_{\pm})$. 

We now compute the second Chern classes $c_2 (\Delta^+ ) $ and 
$c_2 (\Delta^- )$ of the bundles $\Delta^+$ and $\Delta^-$. 
Since we have the relations 
$p_1 (E_+ )  = -4 c_2 (\Delta^+ ) $
and $p_1 (E_- )  = -4 c_2 (\Delta^- ) $
between the first Chern classes and the first Pontryagin classes it is enough 
to compute the first Pontryagin classes of the bundles $E_{\pm}$. 

\proclaim Proposition 3.2.7. Let $E_+ = T\Mg \times_{\rho_+} 
\Lambda_+^2 (\R^4 )$ and 
$E_- = T\Mg \times_{\rho_-} \Lambda_-^2 (\R^4 )$
be the vector bundles induced from the tangent bundle of the manifold $\Mg$ 
by the representations 
$ \rho_{\pm} : SO(4) \rightarrow GL (\Lambda^2_{\pm} (\R ^4 ))$.
The first Pontryagin classes are  
$p_1 (E_+ )= 4(1-2g)$ and  
$ p_1 (E_- ) = - 4(1-2g) $. 

\noindent{\it Proof.} 
$p_1 (E_-)$ follows from $p_1 (E_+)$ by changing the orientation. 
Let $T \subset SO(4)$ be the standard maximal torus 
$T \cong SO(2) \times SO(2)$ in the compact Lie group $SO(4)$  
and let $T' \subset SO(3)$ be the maximal torus 
$$ 
T' = \lbrace \pmatrix {1&0&0\cr 
                 0&\cos \theta & - \sin \theta \cr 
                 0& \sin \theta & \cos \theta \cr } 
\vert \; 0 \leq \theta \leq 2 \pi \rbrace \; . 
$$ 
We compute for the representation $\rho_+$ restricted to the maximal torus 
$T$ inside $SO(4)$ that 
$\rho_+ (T)\subset T'$ and 
$$ 
\rho_+ ( \pmatrix{\cos \theta_1 & - \sin \theta_1 & 0&0 \cr 
                  \sin \theta_1 & \cos \theta_1 & 0& 0 \cr 
                   0& 0& \cos \theta_2 & - \sin \theta_2 \cr 
                   0& 0& \sin \theta_2 & \cos \theta_2 \cr } ) = 
\pmatrix{1&0&0\cr 
         0& \cos (\theta_1 + \theta_2 )& -\sin (\theta_1 + \theta_2 )\cr 
         0& \sin (\theta_1 + \theta_2 ) & \cos (\theta_1 + \theta_2 ) \cr}
\; . 
$$ 
Let $y \in H^1 (T' ; \Z )$ and $x_1 , x_2 \in H^1 (T , \Z )$ be the 
generators. Then the induced map in cohomology is given by 
$$ 
\eqalign{ \rho_+^{\ast} : H^1 (T' ; \Z ) & \rightarrow H^1 (T ; \Z ) \cr 
          \rho_+^{\ast} (y)&= x_1 + x_2 \cr \; . } 
$$ 
By the splitting principle we may assume that 
$T\Mg \cong l_1 \oplus l_2 $ where $l_1 $ and $l_2 $ are $SO(2)$-bundles 
over the manifold $\Mg$ and $x_i = c_1 (l_i ) $ for $i=1,2$. 
By [BoH, Theorem 10.3.b, p. 491] we obtain for the total Pontryagin 
classes of the vector bundle $E_+$ 
$$ 
\eqalignno{ p(E_+ ) = 1 + (x_1 +x_2 )^2 &=  1 + x_1^2 +x_2^2 + 2x_1 x_2 \cr 
 & = 1+p_1 (T\Mg ) + 2 e(T \Mg ) &(3.2.8). \cr }
$$ 
where $p_1 (T \Mg )$ denotes the first Pontryagin class and 
$e (T \Mg )$ the Euler class of the tangent bundle of the manifold 
$\Mg$. 
Hence we obtain from 3.2.8 and lemma 3.2.9 below the equality   
$p_1 (E_+ )  = 2 (2-2(2g))=4(1-2g) $.  
This finishes the proof of proposition 3.2.7. 

\proclaim Lemma 3.2.9. 
(i) $e(T\Mg )= 2-2(2g)$ and 
(ii) $ p_1 (T\Mg )=0$. 

\noindent{\it Proof of lemma 3.2.9.} 
Let $\pi : \Fg \times S(v) \rightarrow \FS=\Mg$ be the canonical projection. 
On the fourth cohomology group $H^4$ with integer coefficients 
the induced map $\pi^{\ast}$ is multiplication by 2. Also we 
obviously obtain for the 
pull back bundle of the tangent bundle $T \Mg$  
the isomorphism $\pi^{\ast} (T \Mg ) \cong T\Fg \times TS(v)$. 
As complex line bundles we write $T\Fg = l_{2g}$ and $TS(v) = l_0$ 
with the dual bundles $\bar{l}_{2g}$ and $\bar{l}_0$. 

For (ii) we obtain: 
$$ 
\eqalign{ 
\pi^{\ast} (p_1(T\Mg ))& = 2p_1 (T \Mg ) \cr 
 & = p_1 (T \Fg \times TS(v) )=-c_2 ((T\Fg \times TS(v) )\otimes \C ) \cr 
 & = -c_2 ((l_{2g}\oplus \bar{l}_{2g})\times (l_0 \oplus \bar{l}_0))\cr  
 & = 1 \times c_1(l_0 )^2 + c_1 (l_{2g})^2 \times 1 \cr 
 & = 0 \quad {\rm evaluated \; on } \quad H^4 (\Fg \times S(v) ; \Z ) \cr}
$$ 
since the classes $c_1 (l_0 )^2$ evaluated on $H^{\ast} (S(v) , \Z )$ and 
$c_1 (l_{2g})^2$ evaluated on $H^{\ast} (\Fg ; \Z )$ are zero. 

For (i) we obtain: 
$$ 
\eqalign{ 
\pi^{\ast} (e(T\Mg ))&=2e(T\Mg ) \cr 
 & =e(T\Fg \times TS(v) )=e(T\Fg )\times e(TS(v) ) \cr 
 & = (2-2(2g))2 \quad {\rm evaluated \; on} \quad H^4 (\Fg \times S(v) ; \Z )
\, . \cr}
$$ 
This gives $e(T \Mg ) = 2-2(2g)$ 
and finishes the proof of lemma 3.2.9. 

\proclaim Corollary 3.2.10. The second Chern classes of the positive and 
negative spinor bundles $\Delta^{\pm}$ over the manifold $\Mg$ are 
given by 
$ c_2 (\Delta^+ ) = 2g -1 $ and 
$ c_2 (\Delta^- ) = 1-2g $.  

We can now set up the equivariant Morse theory. 
Fix an $SU(2)$-invariant metric on the manifold $\Mg$. 

\noindent {\it Remark 3.2.11.} 
Since $H^2 (\Mg ; \Q )=0$ and since the 
bundles $\Delta^{\pm }$ have non zero second Chern classes 
by theorem 3.1 in [FU, p.47] (see also [Pa1, lemma 3.3, p.346]) 
they cannot admit 
connections with stabilizer in the equivariant 
gauge group larger than the constant $\Z_2 \subset SU(2)$ and hence they do not admit 
reducible connections. 

By corollary 3.2.10 if $g=0$ then the bundle $\Delta^+$ does not 
admit anti self dual connections and the bundle $\Delta^-$ does 
not admit self dual connections. If $g \geq 1$ the bundle $\Delta^+$ 
does not admit self dual connections and the bundle $\Delta^-$ does 
not admit anti self dual connections. We fix the following 
conventions: For $g=0$ we study connections on the bundle $\Delta^-$ 
and for $g \geq 1$ we study connections on the bundle $\Delta^+$. 
In the case $g=0$ we define $\B_{SU(2)}$ to be the $SU(2)$-invariant 
connections on the bundle $\Delta^-$ modulo the invariant gauge group. 
In the case $g \geq 1$ we define $\B_{SU(2)}$ to be the $SU(2)$- 
invariant connections on the bundle $\Delta^+$ modulo the invariant 
gauge group. 
Let $\M_{SU(2)}\subset \B_{SU(2)}$ denote the subspace of 
$SU(2)$-invariant anti self dual connections modulo the invariant 
gauge group on the bundle $\Delta^-$ (if $g=0$) or 
$\Delta^+$ (if $g \geq 1$).   
By  changing the orientation we obtain theorems on the existence 
of critical points of the Yang-Mills functional which are not 
self dual on the ``opposite'' bundles. 

In order to use Morse or Lusternik-Schnirelman theory we would like 
to know the weak 
homotopy type of the invariant orbit space 
$ \B_{SU(2)}$. Therefore 
in the next section we will compute the homotopy type of this   
space. 

\proclaim 3.3 Computation of the weak homotopy type of the space  
$\B_{SU(2)}$ .

Let $f_0 \in \Fg$ be a fixed point of the involution 
$\sigma_1 : \Fg \rightarrow \Fg$ defined in 
3.2.1 and $e_0 \in S(v)$ the point $(1,0,0)$ of the two sphere $S(v)$. 
Then on the manifold $\Mg = \FS$ we choose the point 
$[f_0 , e_0 ] \in \FS $ as the base point. 

Since the action of the Lie group $SO(3)$ on the two sphere $S(v)$ 
is transitive the orbit of the $SU(2)$-action on the manifold $\Mg$ through 
the base point $[f_0 , e_0] \in \Mg$ is equal to the subspace 
$\lbrace [f_0 , w ] : w \in S(v) \rbrace \cong \sigma_2 \backslash S(v) \cong 
\R P^2 $. The $SU(2)$-action on the manifold $\Mg$ has 2 types of isotropy 
groups. Every point $[f,w] \in \FS$ is fixed by some circle 
$U(1) \subset SU(2)$. If the point $f \in \Fg$ is not fixed by the involution 
$\sigma_1 : \Fg \rightarrow \Fg$ then the isotropy group of the point 
$[f,w]$ for any $w\in S(v)$ is conjugate to the standard circle 
$$ 
U(1) \cong \lbrace \pmatrix {\exp (i \theta )&0 \cr 
                         0& \exp (-i \theta ) \cr } 
: 0 \leq \theta \leq 2 \pi \rbrace 
\eqno{(3.3.1)} 
$$ 
in the Lie group $SU(2)$. 
If the point $f \in \FS$ is fixed by the involution $\sigma_1 : 
\Fg \rightarrow \Fg$ then the isotropy group of the point $[f,w]$ for 
any $w \in S(v)$ is conjugate to the group $Pin(2) \subset SU(2)$ which is 
generated by the standard circle $U(1) \subset SU(2)$ and the 
element $ j \in Sp(1) \cong SU(2)$.  
The group $Pin(2)$ is the double cover of the group 
$O(2) \subset SO(3)$. 

We now compute the isotropy representations of the two possible isotropy 
groups $U(1)$ and $Pin(2)$ in $SU(2)$ on the fibers of the bundles 
$\Delta^+$ and $\Delta^-$ over the manifold $\Mg$. Recall that the positive 
spinor bundle $\Delta^+$ is a lift of the $SO(3)$-bundle  
$E_+ = P(T\Mg )\times_{\rho_+} \Lambda^2_+ (\R^4 )$. 

Let $m=[f,w]\in \FS =\Mg$ be a point with isotropy group conjugate to the 
standard circle $U(1) \subset SU(2)$ defined in 3.3.1. 
The isotropy representation of the point $m$ on the fiber 
$(T\Mg )_m$ of the tangent bundle $T\Mg$ of the manifold $\Mg$ is 
conjugate to the homomorphism 
$$
\eqalign{\tau : U(1) & \rightarrow SO(4)\cr
\exp (i\theta )&\mapsto \pmatrix{1&0&0&0\cr 
                                 0&1&0&0\cr 
                                 0&0&\cos (2\theta )&-\sin (2\theta )\cr 
                                 0&0&\sin (2\theta )&\cos (2\theta )\cr } 
\qquad . 
\cr}
$$ 
Hence the isotropy representation of the circle $U(1)$ on the fiber $(E_+ )_m$ 
of the bundle $E_+$ is given by the homomorphism 
$$ 
\eqalign{\rho_+ \circ \tau : U(1) &\rightarrow SO(3) \cr 
\exp(i \theta )&\mapsto \pmatrix{1&0&0\cr 
                                 0&\cos (2\theta )&-\sin (2\theta )\cr 
                                 0&\sin (2\theta )&\cos (2 \theta ) \cr} 
\cr} 
\eqno{(3.3.2)} 
$$ 
Since the bundle $\Delta^+$ is the spinor bundle of the bundle $E_+$ for 
some spin structure on the bundle $E_+$ the isotropy representation of the 
point $m$ on the fiber $( \Delta^+ )_m$ 
of the bundle $\Delta^+$ is conjugate to the representation 
$$ 
\eqalign{U(1) & \rightarrow SU(2)\cr 
\exp (i\theta )&\mapsto \pm \pmatrix{\exp (i\theta )&0\cr 
                                      0&\exp(-i\theta )  \cr}  \qquad . 
\cr }
\eqno{(3.3.3) }
$$ 
Since the circle $U(1)$ is connected 
the sign in 3.3.3 has to be constant and equal 
to $+$. Hence the isotropy representation 
of the point $m=[f,w]\in \Mg$ is conjugate to the standard 
inclusion $ U(1) \hookrightarrow SU(2)$ given in 3.3.2.  

Now let $m=[f,w] \in \FS =\Mg$ be a point with isotropy group 
conjugate to the group 
$Pin(2)= \langle j, 
U(1)\rangle \subset 
SU(2)$. The isotropy representation of the group $Pin(2)$ on the fiber 
$(\Delta^+ )_m$ of the positive spin bundle  
is a 2-dimensional complex representation. Since 
restricted to the circle $U(1) \subset Pin(2)$ this representation is 
conjugate to the standard inclusion $U(1) \subset SU(2)$ (by the same argument 
as above) it has to be 
conjugate to the unique irreducible 
2-dimensional complex representation of the group $Pin(2)$ on $\C^2$ 
which is induced from the circle $U(1) \subset Pin(2)$ by the 1-dimensional 
complex representation of weight 1 (or weight $-1$). This 
representation is just given by the standard inclusion $Pin(2) \hookrightarrow 
SU(2)$. The same arguments apply to the isotropy representations of the 
negative spinor bundle. 

Identifying the Lie group $SU(2)$ with the Lie group $Sp(1)$ we can consider 
the bundles $\Delta^+$ and $\Delta^-$ as $SU(2)$-equivariant quaternionic line 
bundles. 
Recall that an $SU(2)$-equivariant quaternionic line bundle 
$\eta = (E \rightarrow \Mg )$ is subordinate to a quaternionic left module 
$V$ if for every point $m \in \Mg$ the isotropy representation 
of the isotropy group $SU(2)_m$ on the fiber $E_m$ of the bundle 
$\eta$ over the point $m \in \Mg$ is contained up to 
isomorphism in the space $V$ viewed as an $SU(2)_m$-module. 
The previous discussion proves the following 
proposition: 

\proclaim Proposition 3.3.4. The $SU(2)$-equivariant quaternionic line bundles 
$\Delta^+$ and $\Delta^-$ over the manifold $\Mg$ are subordinate 
to the non trivial 1-dimensional quaternionic 
$SU(2) \cong Sp(1)$-module $\H$ given by left multiplication of the group 
$Sp(1)$. 

Let $\G^+_{SU(2)}$ denote the $SU(2)$-equivariant gauge group of the bundle 
$\Delta^+$ and let $\G^{+,0}_{SU(2)}$ denote the subgroup of based 
equivariant gauge transformations, i.e. the equivariant gauge transformations on the bundle $\Delta^+$ which are the identity on the fiber 
$(\Delta^+ )_m$ over the base point 
$m=[f_0 ,e_0  ]\in \FS =\Mg$. 
Similarly we define the equivariant gauge 
group $\G^-_{SU(2)}$ and the based equivariant gauge group $\G^{-,0}_{SU(2)}$ 
of the bundle $\Delta^-$. 

\proclaim Lemma 3.3.5. We have isomorphisms 
(i) $\G^+_{SU(2)} / \G^{+,0}_{SU(2)} \cong \Z_2 $
and (ii) $\G^-_{SU(2)} / \G^{-,0}_{SU(2)} \cong \Z_2 $
 
\noindent{\it Proof.} For the proof of (i) 
let $\Delta^+_m$ denote the fiber of the vector bundle 
$\Delta^+$ over the base point $m=[f_0 ,e_0 ]\in \Mg$. Recall that the 
base point $m \in \Mg$ has the isotropy group $Pin(2) \subset SU(2)$ 
since we chose the point $f_0 \in \Mg$ to be a fixed point of the involution 
$\sigma_1 : \Fg \rightarrow \Fg$. 
Define the restriction map 
$$ 
\eqalign{ \Gamma_m : \G^+_{SU(2)} &\rightarrow GL(\Delta^+_m ) \cr 
\phi&\mapsto \phi \vert_m \qquad . \cr }
$$ 
Since the gauge transformation $\phi \in \G^+_{SU(2)}$ is $SU(2)$-equivariant 
it commutes with the action of the isotropy group $Pin(2)$ on the space 
$\Delta^+_m$ and hence $\phi \vert_m$ is a $Pin(2)$-module isomorphism. 
We have seen above that this $Pin(2)$-module 
$\Delta^+_m$ 
is irreducible. Therefore  
by Schur's lemma it must be given by multiplication with some element 
$\lambda \in \C^{\ast}$, i.e. $\phi \vert_m  = \lambda {\rm Id}: \C^2 \rightarrow 
\C^2$. 
But the bundle $\Delta^+$ has structure group $SU(2)$ and by definition the 
gauge transformation $\phi$ is an automorphism of the $SU(2)$-bundle 
$\Delta^+$ and hence $\phi \vert_m \in SU(2)$. This implies $\lambda \in \Z_2$. 
Hence the restriction map $\Gamma_m$ takes values in $\Z_2 \subset SU(2)$ and 
is obviously surjective. Since the kernel of the restriction map $\Gamma_m$ 
is by definition the group $\G^{+,0}_{SU(2)}$ we have proved (i). 
Since the proof of (ii) is the same we have finished the proof 
of lemma 3.3.5. 

Let $\B^{0}_{SU(2)}$ denote the orbit space of 
$SU(2)$-invariant connections on the bundle $\Delta^+$ 
if $g \geq 1$ or $\Delta^-$ if $g=0$  
modulo the action of 
the based invariant gauge group $\G^{+,0}_{SU(2)}$ or 
$\G^{-,0}_{SU(2)}$ respectively.   
Let 
$
\pi : \B^{+,0}_{SU(2)} \rightarrow \B^+_{SU(2)}
$ 
be the canonical projection. Since the bundles 
$\Delta^{\pm}$ do 
not admit a reducible connection (see remark 3.2.11) the map $\pi $ 
is a fibration with fiber 
$( \G^{\pm}_{SU(2)} / \G^{\pm ,0}_{SU(2)})/\Z_2 $. 
By lemma 3.3.5 we have the isomorphism 
$ \G^{\pm}_{SU(2)} / \G^{\pm ,0}_{SU(2)} \cong \Z_2$ and hence the 
projection map $\pi$ 
is in fact a homeomorphism. This proves 

\proclaim Proposition 3.3.6. 
The natural projection map 
$
\pi : \B^{0}_{SU(2)}  \rightarrow \B_{SU(2)} 
$ 
is a homeomorphism. 

In the remainder of this section we compute the weak homotopy type of the based 
invariant orbit space $ \B_{SU(2)}$ using 
theorem 2.3.1. 

\proclaim Theorem 3.3.7. 
There is a weak homotopy equivalence  
$$ 
\B^{0}_{SU(2)}\simeq \prod_{2g} S^1 
$$  

\noindent{\it Proof.} 
Let $\H$ be the 1-dimensional quaternionic representation of 
$SU(2)\cong Sp(1)$ given by left multiplication. By proposition 3.3.4 the 
$SU(2)$-equivariant bundles $\Delta^{\pm}$ are subordinate to this representation. 
Hence theorem 2.3.1 gives a weak homotopy equivalence 
$$ 
\B^{0}_{SU(2)} \simeq \Map^0_{SU(2)} \, (\Mg , G_1 (\H^{\infty}))^{\pm} \; . 
\eqno{(3.3.8)}  
$$ 
Here we have chosen the point $[f_0 ,e_o ]\in \FS$ 
defined before as the base point and a base 
point $\ast \in G_1 (\H^{\infty})$ 
according to theorem 2.2.2. Both points are fixed by 
the group $Pin(2) \subset SU(2)$. The superscript $\pm$ 
on the right hand side of 3.3.8 denotes 
the components of maps which classify the $SU(2)$-equivariant positive 
spinor bundles $\Delta^{\pm}$ over the manifold $\Mg = \FS$. 

One checks that any line $\ast =[(x_0 ,x_1 , \dots )] $ in the infinite 
dimensional Grassmannian $G_1 (\H^{\infty} )$ which is fixed by the 
group $Pin(2) \subset SU(2)$ is already fixed by the whole group $SU(2)$ and 
hence is a fixed point of the $SU(2)$-action on the space 
$G_1 (\H^{\infty})$. Hence any map $f:\Mg \rightarrow G_1 (\H^{\infty})$ 
which is base point preserving maps the whole orbit through the base point 
$[f_0 ,e_0 ]$ to the base point $\ast \in  G_1 (\H^{\infty})$, i.e. 
$f([f_0 ,x])=\ast$ for any $x \in S(v)$. This implies the homeomorphisms  
$$ 
\Map^0_{SU(2)} \, (\FS,\, G_1(\H^\infty))\cong \Map^0_{SU(2)}\, 
({\FS\over \{f_0\}\times_\sigma S(v)}, \, G_1(\H^\infty)) \; . 
\eqno{(3.3.9)} 
$$

Let $S(v)_+ =S(v) \cup \{ + \}$ denote the two sphere $S(v)$ with an additional 
disjoint base point $\{ +\}$. 
We choose the point $f_0 \in \Fg$ as the base point in the Riemann surface 
$\Fg$ 
and define the smash product 
$$ 
\Fg \wedge S(v)_+ = {\Fg \times S(v)_+ \over \Fg \times \{ +\} \cup 
\{ f_0 \} \times S(v)_+ }\qquad . 
$$ 
We extend the involution $\sigma_2 : S(v) \rightarrow S(v)$ given by the 
reflection through the origin 
to the space $S(v)_+$ by letting it act trivially on 
the additional base point. 
The involution $\sigma : \Fg \times S(v) \rightarrow \Fg \times S(v)$ given by 
$\sigma (v,w)= (\sigma_1 (v) , \sigma_2 (w))$ extends to the space 
$\Fg \times S(v)_+$ and it preserves the subspace 
$  \Fg \times \{ +\} \cup 
\{ f_0 \} \times S(v)_+ \subset \Fg \times S(v)_+$. Hence we obtain an induced 
involution $\sigma : \Fg \wedge S(v)_+ \rightarrow \Fg \wedge S(v)_+$. 
Let $\Fg \wedge_{\sigma} S(v)_+$ denote the quotient.

We extend the $SU(2)$-action on the two sphere $S(v)$ 
to the space $S(v)_+$ by letting it act trivially 
on the additional base point $+$. This induces an action on the smash 
$\Fg \wedge S(v)_+ $ which commutes with the action of the involution 
$\sigma$. Hence it induces an $SU(2)$-action on the quotient 
 $\Fg \wedge_{\sigma} S(v)_+$. 
The canonical mapping 
$$ 
\eqalign{  
\varepsilon : \quad {\FS\over \{f_0\}\times_\sigma S(v)} \; 
& \longrightarrow \; \Fg \wedge_{\sigma} S(v)_+ \cr 
\varepsilon ([x,y])&= [x,y]\cr}
$$ 
is equivariant with respect to the $SU(2)$-action and induces an 
$SU(2)$-equivariant homeomorphism. 
Together with 3.3.9 this gives a homeomorphism 
$$ 
\Map^0_{SU(2)} \, (\FS;\, G_1(\H^\infty)) \cong 
\Map^0_{SU(2)} \, (\Fg \wedge_{\sigma} S(v)_+ \, , G_1 (\H^{\infty})) \quad . 
\eqno{} 
$$ 

A map $f: \Fg \wedge_\sigma S(v)_+ \rightarrow G_1 (\H^{\infty})$ 
is given by a map 
${\tilde f} : \Fg \wedge S(v)_+ \rightarrow G_1 (\H^{\infty})$ which is 
invariant under the action by the involution $\sigma$. (The  
space $ G_1 (\H^{\infty})$ is given the trivial $\sigma$-action). 
The adjoint map 
${\tilde f}^{{\rm ad}} : \Fg \rightarrow \Map^0\, (S(v)_+ \, ,G_1 ( \H^{\infty}))$ 
is then invariant under a $\Z_2$-action where the group $\Z_2$ 
acts on the Riemann surface by the action of the involution $\sigma_1$ 
and on the mapping space $ \Map^0\, (S(v)_+ \, ,G_1 ( \H^{\infty}))$ by the 
induced action on maps given by the $\sigma_1$-action on the space $S(v)_+$ 
and the trivial action on the Grassmannian $G_1 ( \H^{\infty})$. 
Since the $SU(2)$-and $\sigma_2$-actions on the space $S(v)_+$ commute we 
obtain the homeomorphism 
$$ 
\Map^0_{SU(2)}\, (\FS \, , G_1 (\H^{\infty}) \cong 
\Map^0_{\Z_2}\, (\Fg\, ,\Map^0_{SU(2)}\, (S(v)_+ \, ,G_1 ( \H^{\infty}))) \; . 
\eqno{(3.3.10)}
$$

We now analyze the mapping space 
$\Map^0_{SU(2)}\, (S(v)_+ \, ,G_1 ( \H^{\infty}))$ together with the 
above $\Z_2$-action. 
We first need to prove 

\proclaim Lemma 3.3.11. 
There is a homeomorphism 
$ \g {1}^{U(1)} \cong G_1 (\C ^{\infty}) $.

{\it Proof of lemma 3.3.11.} 
Let $\tau : G_1 (\C^ {\infty}) \hookrightarrow \g {1}$ be the map induced 
by the natural inclusion $ \C \hookrightarrow \H$. 
Recall that we view the space $\H^{\infty}$ as a right quaternionic 
vector space endowed with the $U(1)$-action given by left 
multiplication of weight p. Hence the image 
$ \tau ( G_1 (\C^ {\infty}))$ lies in the fixed point set 
$ \g {1}^{U(1)}$ of the induced circle action on 
the Grassmannian $\g {1}$ and the map 
$\tau$ induces a continuous map $\tilde{\tau}: G_1 (\C^{\infty}) 
\hookrightarrow \g {1}^{U(1)}$. 
Let $w= [x_0 , x_1 , \ldots ]$ be a quaternionic line in $\H^{\infty}$ 
fixed under the circle action. Without loss of generality we may assume 
$x_0 =1$. Given $\lambda \in U(1)$ there is an element 
$\alpha (\lambda ) \in Sp(1)$ such that 
$\lambda^p x_i = x_i \alpha (\lambda )$ for all $i=0,1, \ldots$. For 
$i=0$ this gives $\alpha (\lambda )= \lambda^p $. Hence 
$\lambda^p  x_i = x_i \lambda^p$ for all $i=0,1, \ldots$ and hence 
$x_i$ lies in the centralizer of $U(1)$ in the quaternions $\H$ 
which is equal to $\C$. Hence $x_i \in \C$ for all $i=0,1, \ldots$ 
and the element $w= [x_0 , x_1 , \ldots ]$ lies in the image of the 
map $\tilde{\tau} : G_1 (\C ^{\infty})\rightarrow \g {1}^{U(1)}$. The  
map $\tilde{\tau} $ induces a homeomorphism $G_1 (\C ^{\infty}) \cong 
\g {1} ^{U(1)}$. 
This finishes the proof of lemma 3.3.11. 

\proclaim Lemma 3.3.12. 
There is a $\Z_2$-equivariant homeomorphism 
$$
\Map^0_{SU(2)}\, (S(v)_+ \, ,G_1 ( \H^{\infty})) \cong G_1 ( \C^{\infty}) 
$$  
where the group $\Z_2$ acts on the infinite dimensional 
Grassmannian $G_1 (\C^{\infty})= \C P^{\infty}$ by complex conjugation on coordinates. 

\noindent{\it Proof.}  
The inclusion $S(v) \hookrightarrow S(v)_+$ induces a homeomorphism 
$$ 
\Map^0_{SU(2)}\, (S(v)_+ \, ,G_1 ( \H^{\infty})) \cong 
\Map_{SU(2)} (S(v) \, , G_1 ( \H^{\infty})) \; . 
$$ 
Let $N \in S(v)$ denote the north pole. The evaluation map 
$$ 
\eqalign{ \Map_{SU(2)} (S(v) \, ,G_1 ( \H^{\infty})) \; & \rightarrow 
G_1 ( \H^{\infty})^{U(1)} \cr 
f&\mapsto f(N) \cr }
\eqno{(3.3.13)} 
$$ 
induces an $SU(2)$-equivariant homeomorphism. 
Here we identify the two sphere $S(v)$ with the homogeneous 
space $SU(2)/U(1)$. 
By lemma 3.3.11 there is a homeomorphism 
$G_1 (\C^\infty ) \cong  G_1 ( \H^{\infty})^{U(1)}$. 

We now compute the 
$\Z_2$-action on the infinite dimensional complex Grassmannian 
$G_1 (\C^\infty )$ induced by the homeomorphism in 3.3.13. 
Let $\tau$ be a generator of $\Z_2$. For any map 
$f \in  \Map_{SU(2)} (S(v) \, ,G_1 ( \H^{\infty}))$ the map 
$\tau \cdot f$ is given by $\tau \cdot f (w)= f( \sigma_2 (w))$. 
Hence $\tau \cdot f (N)= f(\sigma_2 (N)) = f(-N) =f(S)$. 
Here we denote by $S \in S(v)$ 
the south pole. Let $A \in SU(2)$ be an element such that 
$A(N)=S$. 
Since the map $f$ is equivariant with 
respect to the $SU(2)$-action we obtain 
$\tau \cdot f (N)= f( \sigma_2 (N))= f(S) = f(A(N)) = A(f(N))$. 
If $B$ is another element in $SU(2)$ such that $B(N)=S$ then 
$B=C \cdot A$ for some element $C \in U(1) \subset SU(2)$. 
Here $U(1) \subset SU(2)$ is the stabilizer of the $SU(2)$-action on $S(v)$ at 
the north pole $N$. 
But since 
$f(N) $ lies in the fixed point set of the $U(1)$-action on the infinite 
dimensional Grassmannian $G_1 (\H^\infty )$ we obtain 
$B(f(N))=A(f(N))$. Hence the $\Z_2$-action on the space 
$G_1 ( \H^{\infty}))^{U(1)}$ is given by left multiplication by any element 
$A \in SU(2)$ such that $A(N)=S$. 

We now identify $SU(2)$ with $Sp(1)$ and the two sphere $S(v)$ with the 
purely imaginary quaternions of norm 1. Then the 
$SU(2)$-action on the 2-sphere $S(v)$ is given by conjugation on 
quaternions. 
The north pole $N$ is identified with the quaternion $k $. 
Since $jk {\bar j}= -jkj=-k$ we can choose the element $A \in SU(2)$ 
to be conjugation by the quaternion $j \in Sp(1)$. Hence the $\Z_2$-action 
on the space $G_1 (\H^\infty )^{U(1)}$ induced by the homeomorphism 3.3.13  
is given by left multiplication by 
the quaternion $j$. Let $[x_0 ,x_1 , \dots ]$ be a line in the Grassmannian 
$G_1 (\H^\infty )$ invariant under the $U(1)$-action. By lemma 3.3.11
we may assume 
that the coordinates $x_i \in \H$ are all complex, i.e. 
$x_i \in \C$ for all $i$. Hence 
$$
j \cdot [x_0 ,x_1 , \dots ] = [jx_0 ,jx_1 , \dots ]
= [\bar{x}_0 j, \bar{x}_1 j , \dots ] =  [\bar{x}_0 ,\bar{x}_1 , \dots ] 
$$ 
and therefore the $\Z_2$-action on the space 
$G_1 (\H^\infty )^{U(1)} \cong G_1 (\C^\infty)$ is given by complex 
conjugation. This finishes the proof of lemma 3.3.12. 

Hence we have proved 

\proclaim Proposition 3.3.14. 
There is a canonical homeomorphism 
$$ 
\Map^0_{SU(2)}\, (\FS , \, G_1 (\H^\infty )) \cong 
\Map^0_{\Z_2} \, (\Fg , G_1 (\C^\infty )) 
$$ 
where the $\Z_2$-action on the Riemann surface $\Fg$ is given by the 
involution $\sigma_1$ defined in 3.2.1 and on the infinite dimensional 
Grassmannian $ G_1 (\C^\infty )$ by complex conjugation. 

We now compute the homotopy type of the 
$\Z_2$-equivariant mapping space 
$ \Map^0_{\Z_2} \, (\Fg , G_1 (\C^\infty )) $. 
The inclusion of the ``gluing circle''  $S^1$ into the 
Riemann surface 
$\Fg = ({\cal F}_g -D^2 )\cup_{S^1} ({\cal F}_g -D^2 ) $ gives a 
$\Z_2$-equivariant cofiber sequence 
$$ 
S^1 \hookrightarrow \Fg \rightarrow \Z^+_2 \wedge {\cal F}_g 
\eqno{(3.3.15)} 
$$ 
where the non trivial element $-1 \in \Z_2$ acts on the surface 
$\Fg$ by reflection on the symmetry plane. 
Applying the equivariant mapping functor 
$\Map^0_{\Z_2}\, (\underline{\quad}, G_1 (\C^\infty ))$ to 3.3.15 
gives a fibration 
$$ 
\Map^0_{\Z_2}\, (S^1 , \, G_1 (\C^\infty )) \leftarrow 
\Map^0_{\Z_2}\, (\Fg, \, G_1 (\C^\infty )) \leftarrow 
\Map^0_{\Z_2}\, (\Z_2^+ \wedge {\cal F}_g , G_1 (\C^\infty )) 
\; . 
\eqno{(3.3.16)} 
$$ 

Since the group $\Z_2$ acts trivially on the circle $S^1$  
the base space of the fibration 3.3.16 is homeomorphic to the space 
$\Omega ( G_1 (\C^\infty )^{\Z_2})$. 
Similarly as in lemma 3.3.11 one computes that the inclusion 
$ G_1 (\R^\infty ) \hookrightarrow G_1 (\C^\infty )$ 
induced by the canonical inclusion 
$\R \hookrightarrow \C$ 
induces a homeomorphism 
$ (G_1 (\C^\infty ))^{\Z_2} \cong G_1 (\R^\infty )$. 
Hence the base space of fibration 3.3.16 is homotopy equivalent to the 
2-point space $\Z_2$, 
the total space is homotopy equivalent 
to the disjoint union of two copies of the 
fiber and we obtain the homotopy equivalence  
$$ 
\Map^0_{\Z_2}\, (\Fg , G_1 (\C^\infty )) \simeq  
\bigsqcup_{2 \; {\rm copies}} 
\Map^0_{\Z_2}\, (\Z_2^+ \wedge {\cal F}_g , G_1 (\C^\infty )) \; . 
\eqno{(3.3.17)} 
$$
There is a natural homeomorphism 
$$
\Map^0_{\Z_2}\, (\Z_2^+ \wedge {\cal F}_g , G_1 (\C^\infty ))
\cong 
\Map^0 \, ({\cal F}_g ,  G_1 (\C^\infty )) \; . 
$$ 
Since $G_1 (\C^\infty ) =BU(1)$ 
there is a canonical homotopy equivalence  
$$
\Map^0 \, ({\cal F}_g ,  G_1 (\C^\infty ))
\simeq \Z \times \prod_{2g} S^1 \; . 
$$ 
Here we have switched from the notation $U(1)$ for the circle as a group 
to $S^1$ for the circle as a topological space since a mapping space is not 
necessarily a group. 
Therefore we obtain the homotopy equivalence 
$$ 
\Map^0_{\Z_2}\, (\Fg , G_1 (\C^\infty )) \simeq 
\bigsqcup_{2 \; {\rm copies}} (\Z \times \prod_{2g} S^1 ) \; . 
\eqno{(3.3.18)} 
$$ 
Hence, using proposition 3.3.14, we obtain 
a homotopy equivalence 
$$
\Map^0_{SU(2)} \, (\FS , G_1 (\H^\infty )) \simeq  
\bigsqcup_{2 \; {\rm copies}} (\Z \times \prod_{2g} S^1 ) 
\eqno{(3.3.19)}  
$$
and each component of the mapping space in 3.3.19 is homotopy equivalent 
to the $2g$-fold product of the circle $S^1$. Together with 3.3.8 
this finishes the proof of theorem 3.3.7.

\proclaim 3.4 Proof of the main theorem. 

Recall that we defined the space $\M_{SU(2)} \subset
\B_{SU(2)}$ to be the space of invariant anti self dual 
connections modulo the invariant gauge group on the bundle 
$\Delta^-$ if $g=0$ and on the bundle $\Delta^+$ if  
$g \geq 1$. 
Since by remark 3.2.11 the bundles $\Delta^+$ and $\Delta^-$ do not admit 
reducible connections the invariant moduli space $\M_{SU(2)}$ 
consists of irreducible connections. By proposition 3.1 (p.446) in [Ba] 
(see also theorem 4.6 (p.248) and section 5 in [Cho])  
for a set of 
$SU(2)$-invariant $C^q$-metrics ($q \geq 1$) on the manifold $\Mg$, 
open and dense in the set of all invariant metrics,  the 
invariant moduli space of anti self dual connections $\M_{SU(2)}$
is a smooth (possibly empty) manifold of a dimension 
which can be computed using the Atiyah-Segal-Singer fixed point formula. 

\noindent{\it Remark 3.4.1.} Proposition 4.1 in [Ba] is stated for the case $H=U(1)$ 
but the proof carries over word for word to the case of any compact Lie group. 

\noindent{\it Remark 3.4.2.} 
The following argument (given to me by J. R\aa de) shows that 
Uhlenbeck's generic metrics theorem 
([FU, Theorem, 3.17, p.59]) and its equivariant versions in 
[Ba] and [Cho] are also true in our Morrey space completions. 
Any anti self dual connection $A \in \M_{SU(2)}$ (completed in the 
$L^{1,2}_{3}$-norm) is smooth since it is a Yang-Mills connection. So we 
get the same moduli spaces no matter whether we work with Morrey or 
Sobolev spaces. In a neighborhood of an anti self dual connection $A$ 
the moduli space $\M_{SU(2)}$ is a manifold if and only if the operator 
$ D_A^+ : \Omega ^2_+ ({\rm ad} (\eta )) \rightarrow 
\Omega ^1 ({\rm ad} (\eta )) $ has trivial null space. 
Here  
${\rm ad} (\eta )$ is the vector bundle associated to the given $SU(2)$-bundle 
$\eta$ via the adjoint representation of $SU(2)$. 
$D_A^+$ is an elliptic operator with smooth coefficients so 
anything in the null space is smooth. Hence the condition 
that the operator $D_A^+$ has trivial null space 
is independent 
of the choice of the function spaces.  

We now compute the formal dimension of the moduli space of invariant 
anti self dual connections. 
Assume that the moduli space $\M_{SU(2)}$ is not empty.  
By standard gauge theory and an argument in [AHS, p.444 and 445] it follows 
that 
for a generic invariant metric 
the equivariant index of the Dirac operator 
$$ 
{D \! \! \! \! /}^-_E : \Gamma (\Delta^- \otimes \Delta^+ 
\otimes ({\rm ad}(\eta )\otimes \C )) 
\rightarrow 
 \Gamma (\Delta^+ \otimes \Delta^+ \otimes ({\rm ad}(\eta )\otimes \C )) 
\eqno{(3.4.3)} 
$$ 
is an actual representation. 
Here $E= \Delta^+ \otimes ({\rm ad}(\eta )
\otimes \C )$ denotes the coefficient 
bundle. 
Also by standard equivariant gauge theory 
for a generic invariant metric 
the dimension of the 
moduli space $\M_{SU(2)}$ is equal to the dimension of the   
trivial representation contained in the 
representation $ {\rm ind} ({D \! \! \! \! /}^-_E )$. 

Let 
$$
h= \pmatrix{e^{i\theta}&0\cr 
            0&e^{-i\theta}\cr}
$$ 
be an element in the standard 
maximal torus $U(1) \subset SU(2)$. If the angle $\theta$ is irrational then the 
closure of the cyclic group generated by the element $h$ in $SU(2)$ 
is just the torus $U(1)$ defined in 3.3.2. The action of the element $h$ on the 
two sphere $S(v)$ fixes the point $e_0^+ =(1,0,0)$ and 
$e_0^- =(-1,0,0)$. Let 
$\Mg^h$ denote the fixed point set of the element 
$h \in SU(2)$ acting on the manifold $\Mg$. Then 
$\Mg^h =\Fg \times_{\sigma} \{ e_0^+ , e_0^- \}\cong \Fg$ where the 
isomorphism is induced by the inclusion 
$\Fg \hookrightarrow \Mg , \; f\mapsto [f,e_0^+ ]$. 

The tangent bundle $T\Mg$ of the manifold $\Mg$ restricted to the fixed 
point set $\Mg^h$ splits 
$U(1)={\overline {\langle h \rangle}}$-equivariantly as 
$$ 
T\Mg \cong T\Fg^0 \oplus N(2\theta ) 
\eqno{(3.4.4)} 
$$ 
where $T\Fg^0$ denotes the real tangent bundle of the surface $\Fg$ 
together with the trivial $U(1)$-action and $N(2\theta )$ denotes the real two 
dimensional trivial bundle. Here the generator $h=e^{i\theta}$ of the 
circle $U(1)$ acts on the bundle $N(2\theta )$ by the matrix 
$$ 
\pmatrix{\cos (2\theta )& -\sin (2\theta )\cr 
          \sin(2\theta )& \sin(2\theta )\cr } \qquad . 
$$ 
The bundle $N(2 \theta )$ is the equivariant normal bundle of the fixed 
point set $\Mg^h \cong \Fg$ inside the 4-manifold $\Mg$. Both bundles in 3.4.4 
are the underlying real bundles of complex line bundles on the surface 
$\Fg$. 

By the $G$-Index Theorem ([AS, 5.4, p.572] 
we obtain 
$$ 
{\rm ind}_h ({\not  \! \! D}^-_E )= \pm {\rm ch}_h(\Delta^+ ) \, 
{\rm ch}_h ({\rm ad}(\eta )\otimes \C ) \, 
(- {\hat A}_{2\theta }(N(2\theta )))\, {\hat A}(\Fg )\; [\Fg] \quad . 
\eqno{(3.4.5)} 
$$ 
Here ${\rm ind}_h ({\not \! \! D}^-_E )$ denotes the index of the operator 
${\not \! \! D}^-_E$ evaluated on the element 
$h \in U(1) \subset SU(2)$, ${\rm ch}_h$ the 
equivariant Chern character, ${\hat A}(\Fg)$ the ${\hat A}$-genus of the 
surface $\Fg$ and $ {\hat A}_{2\theta }(N(2\theta )))$ a certain 
characteristic class to be computed later. 

Let $x\in H^2 (\Fg, \Q )$ be the first Chern class of the complex tangent 
bundle of the Riemann surface $\Fg$. We compute: 

$$
\eqalign{ 
{\hat A}(\Fg )&={x\over e^{x\over 2}-e^{-x\over 2}} \cr 
{\hat A}(N(2\theta ))&= {e^{{1\over 2}(2i\theta )} \over 
e^{2i\theta}-1} = 
{1\over e^{i\theta}-e^{-i\theta }} \cr  
{\rm ch}_h (\Delta^{\pm} ) &= 
e^{-x\over 2} e^{\mp i\theta }+e^{x\over 2} e^{{\pm}i\theta } \cr   
{\rm ch}_h ({\rm ad}(\Delta^{\pm} )\otimes \C )&= 
e^{-x}e^{{\mp}2i\theta}+1+e^{x}e^{{\pm}2i\theta} \cr  }
$$ 
 
Hence for $g=0$ and $\eta =\Delta^- $ we obtain 
$$ 
\eqalign{ 
{\rm ind}_h ({\not \! \! D}^-_E )& 
=- {(e^{-{x\over 2}}e^{-i\theta }+e^{{x\over 2}} e^{i\theta})
(e^{-x}e^{2i\theta }+1+e^{x}e^{-2i\theta })x \over 
(e^{i\theta }-e^{-i\theta })(e^{{x\over 2}}-e^{-{x\over 2}})}
\quad [S^2 ] \cr 
&= e^{2i\theta }+3+e^{-2i\theta } \qquad . \cr }  
\eqno{(3.4.6)} 
$$ 
For $g\geq 1$ and $\Delta^+$ we obtain 
$$ 
\eqalign{ 
{\rm ind}_h ({\not \! \! D}^-_E )& 
=- {(e^{-{x\over 2}}e^{-i\theta }+e^{{x\over 2}} e^{i\theta})
(e^{-x}e^{-2i\theta }+1+e^{x}e^{2i\theta })x \over 
(e^{i\theta }-e^{-i\theta })(e^{{x\over 2}}-e^{-{x\over 2}})}
\quad [\Fg ] \cr 
&= (2g-1)(3e^{2i\theta }+5+3e^{-2i\theta } )\qquad . \cr }  
\eqno{(3.4.7)} 
$$

\noindent{\it Remark 3.4.8.} Careful investigation shows that 
the sign in formula 3.4.5 turns out to be $+1$. 

Hence 
the dimension of the trivial representation $V_{2g} \subset 
{\rm ind}({\not \! \! D}^-_E )$ is   
$$ 
\dim (V_{2g})=\cases{2,&if $g=0$\cr 
                     2(2g-1),&if $g\geq 1 \; .$ \cr} 
$$ 
This gives 

\proclaim Theorem 3.4.9. 
If $g=0$ then either the moduli space of invariant anti self dual
connections $\M_{SU(2)}$ of the bundle $\Delta^-$ is empty or, for a
generic $SU(2)$-invariant metric a manifold of dimension 2.  If $g\geq
1$ then either the moduli space of invariant anti self dual
connections $\M_{SU(2)}$ on the bundle $\Delta^+$ is empty or, for a
generic $SU(2)$-invariant metric a manifold of dimension $2(2g-1)$.

We are now able to prove the main theorem of this section: 

\proclaim Theorem 3.4.10. 
Fix $g=0,1,2,\dots$ and let $\Delta^+$ and $\Delta^-$ be 
the positive and negative $SU(2)$-equivariant spinor bundles over the 
4-manifold $\Mg =\FS$ of second Chern class 
$c_2 (\Delta^+ )=2g-1 $ and $c_2(\Delta^- )=1-2g$. 
Fix an $SU(2)$-invariant metric. 
\item{(i)} The Yang-Mills functional restricted to the invariant orbit space 
${\cal YM}: \B_{SU(2)}\rightarrow \R$
has at least $2g+1$ critical points on each of 
the bundles $\Delta^+$ and $\Delta^-$. 
\item{(ii)} In the case $g=0$ for a generic $SU(2)$-invariant metric 
the 
critical point on the bundle $\Delta^+$ cannot be self dual and on the bundle 
$\Delta^-$ it cannot be anti self dual. 
In the case $g\geq 2$ for a generic $SU(2)$-invariant metric 
at least one of 
the critical points on the bundle $\Delta^+$ cannot be anti self dual and at 
least one of the critical points on the bundle $\Delta^-$ cannot be self dual.

By Palais' principle of symmetric criticality in [Pal2] 
each of the critical points in (i) is a critical 
point of the Yang-Mills functional on the non-equivariant orbit space 
${\cal YM}:\B^+ \rightarrow \R$,  
i.e. they are Yang-Mills connections.

\noindent{\it Proof.} 
Since the Yang-Mills functional is bounded from below and 
the invariant orbit space $\B_{SU(2)}$ 
is a complete Riemannian manifold 
we can apply the classical theorem on Lusternik-Schnirelman theory by 
Palais ([Pal1, Theorem 7.2, p.131]). 
For a topological space $X$ let 
${\rm cat}(X)$ be the category of $X$, that is  
the least integer $n$ so that X can be covered by $n$ 
closed contractible subsets of $X$. Also define ${\rm cl}(X)$, the 
cuplength of the space $X$, to be the largest integer $k$ such that there 
exists a ring $R$ and cohomology classes 
$\alpha_1 , \dots , \alpha_{k-1} \in H^{\ast}(X,R) $ with 
positive dimension such that their cup product does not vanish. 
One always has the inequality 
${\rm cat}(X) \geq {\rm cl}(X) $.  

Palais' theorem says that the Yang-Mills functional on the space 
$\B_{SU(2)}$ has at least ${\rm cat}(\B_{SU(2)})$ 
critical points. 
Theorem 3.3.7, proposition 3.3.6 
and the inequality 
${\rm cat}(X) \geq {\rm cl}(X) $ 
imply that the Yang-Mills functional on the 
space $\B_{SU(2)}$ has at least 
${\rm cl} (\prod_{2g} S^1 )=2g+1$ critical points. 
This proves assertion (i). 

By theorem 3.4.9 the invariant moduli space 
$\M_{SU(2)}$ is	either empty or a compact closed 
manifold of the dimension given in theorem 3.4.9.    
If this manifold  is empty 
then there is nothing to prove. If it is 
not empty we proceed as follows: 

In the case $g=0$ the moduli space $\M_{SU(2)}$ of anti self dual 
$SU(2)$-invariant connections on the bundle $\Delta^-$ has dimension 2. Hence 
$H^2 (\M_{SU(2)}; \Z_2 )\cong \Z_2$. But $H^2 (\B_{SU(2)}; \Z_2 ) =0$ and 
hence the inclusion $\M_{SU(2)}\hookrightarrow \B_{SU(2)}$ cannot be a 
homotopy equivalence. Since the Palais-Smale 
condition is satisfied 
the critical point from (i) cannot lie in the subspace $\M_{SU(2)}$ and 
hence it cannot be anti self dual. 

The case $g \geq 2$ uses the same argument: 
Since $\dim (\M_{SU(2)})=2(2g-1)$ we obtain 
$H^{2(2g-1)}(\M_{SU(2)}; \Z_2 )\cong \Z_2 $. 
Since $\B_{SU(2)}\simeq \prod_{2g} S^1 $ and 
$2(2g-1) > 2g$ we obtain 
$H^{2(2g-1)}(\B_{SU(2)}; \Z_2 )=0$. 
Hence the inclusion $\M_{SU(2)} \hookrightarrow \B_{SU(2)}$ cannot be a 
homotopy equivalence and the Palais-Smale condition guarantees the existence 
of a critical point which does not lie in the subspace 
$\M_{SU(2)}$. 
This proves the case 
$g\geq 2$. The other cases are proved by changing the orientation. This 
finishes the proof of assertion (ii). 

\bigskip 

\proclaim 5~ References. 

\item{[AB]} Atiyah, M.F., Bott, R.: The Yang-Mills equation over 
Riemann surfaces, Phil. Trans. R. Soc. Lond. A 308 (1982), 523-615. 

\item{[AHS]} Atiyah, M.F., Hitchin, N.J., Singer, I.M.: 
Self-Duality in four dimensional Riemannian geometry, 
Proc. R. Soc. Lond. A 362 (1978), 425-61. 

\item{[AS]} Atiyah, M.F., Singer, I.M.: The index of elliptic 
operators: III, Ann. of Math. 87 (1968), 546-604.  

\item{[BoH]} Borel, A., Hirzebruch, F.: Characteristic classes and 
homogeneous spaces I, Amer. J. Math. 80 (1958), 458-538. 

\item{[Ba]} Braam, P.J.: Magnetic Monopoles on three-manifolds, J. 
Differential Geometry Vol. 30 (1989), 425-464.  

\item{[Cho]} Cho, Y. S.: Finite group actions on the moduli space of self 
dual connections I, Trans. AMS, Vol. 323, No. 1 (1991), 233-261. 

\item{[DK]}{Donaldson, S.K., Kronheimer, P.B.: ``The Geometry of 
Four-Manifolds'', 
Oxford mathematical monographs, 1990, Oxford: Oxford University Press.} 

\item{[FU]} Freed, D.S., Uhlenbeck, K.K.: ``Instantons and Four-Manifolds'', 
MSRI publications Vol. 1, 1984, New York-Berlin-Heidelberg: 
Springer Verlag. 

\item{[LM]} Lawson, H.B. Jr., Michelsohn, M.-L.: ``Spin Geometry'', 
Princeton Mathematical Series Vol. 38, 1989, Princeton: Princeton 
University Press. 

\item{[MS]} Milnor, J.W., Stasheff, J.D.: ``Characteristic Classes'', 
Annals of Mathematical Studies No. 76, 1974, Princeton: 
Princeton University Press.   

\item{[Pal1]} Palais, R.S.: Lusternik-Schnirelman Theory on 
Banach Manifolds, Topology Vol. 5 (1966), 115-132. 

\item{[Pal2]} Palais, R.S.: The Principle of Symmetric Criticality, 
Comm. Math. Phys. 69 (1979), 19-30. 

\item{[Pa1]} Parker, T.H.: A Morse theory for equivariant Yang-Mills, 
Duke Mathematical Journal, Vol. 66, No. 2 (1992), 337-356. 

\item{[Pa2]} Parker, T.H.: Non-minimal Yang-Mills fields and 
dynamics, Invent. math. 107 (1992), 397-420. 

\item{[R\aa ]} R\aa de, J.: Compactness Theorems for Invariant Connections,  
submitted for publication.

\item{[Se]}{Segal, G.B.: Equivariant K-theory, Publ. Math. Inst. 
Hautes \' Etudes Sci. 34 
(1968), 129-151.} 

\item{[SS]} Sadun, L., Segert, J.: Non-self-dual Yang-Mills connections 
with quadrupole symmetry, Comm. Math. Phys. 145 (1992), 363-391.  

\item{[SSU]} Sibner, L.M., Sibner, R.J., Uhlenbeck, K.K. : 
Solutions to Yang-Mills equations that are not self-dual, 
Proc. Natl. Acad. Sci. USA Vol. 86 (1989), 8610-8613.  

\item{[Wa]}{Wasserman, A.G.: Equivariant differential topology, Topology 
Vol. 8, 127-150.}  

\item{[Wan]} {Wang, H.-Y.: The existence of nonminimal solutions to the 
Yang-Mills equation with group $SU(2)$ on $S^2 \times S^2$ and 
$S^1 \times S^3$, J. Differential Geometry Vol. 34 (1991), 701-767.}

\bigskip 

\noindent U. Gritsch, DPMMS, University of Cambridge, 16 Mill Lane, 
Cambridge, CB2 1SB, U.K.. 

\noindent e-mail: {\tt U.Gritsch@dpmms.cam.ac.uk}

\end